\documentclass[12pt]{article}
\usepackage[utf8]{inputenc}
\usepackage[margin=2cm]{geometry}
\usepackage{amsmath, amssymb}
\usepackage{graphicx}
\usepackage{xcolor}
\usepackage{comment}
\usepackage[authoryear]{natbib}
\usepackage{multirow}
\usepackage{multicol}
\usepackage{longtable}
\usepackage{booktabs}

\newcommand{\x}{\mathbf{x}}

\usepackage{hyperref}

\usepackage{titlesec}
\titleformat*{\section}{\normalsize\bfseries}% solves Package titlesec Error: Entered in horizontal mode
\titleformat{\paragraph}
{\normalfont\normalsize\bfseries}{\theparagraph}{1em}{}
\titlespacing*{\paragraph}
{0pt}{3.25ex plus 1ex minus .2ex}{1.5ex plus .2ex}
\usepackage{tikz}
\usetikzlibrary{shapes, arrows}
\usepackage{booktabs}%for excel to latex convention
\usepackage{graphicx}%for \resizebox
\usepackage{tabularx}
\usepackage{mathtools}
\usepackage{hyperref}

\usepackage[normalem]{ulem} %for using \sout{}

\usepackage{authblk}

\title{Decision support for sustainable forest harvest planning using multi-scenario multiobjective robust optimization} 

\author[1]{Babooshka Shavazipour}%[orcid=0000-0002-6516-4423]
\author[2]{Lovisa Engberg Sundström}
\affil[1]{University of Jyvaskyla, Faculty of Information Technology, P.O. BOX 35 (Agora), FI-40014 University of Jyvaskyla, Finland}
\affil[2]{Skogforsk, Uppsala, Sweden}

\begin{document}
\let\WriteBookmarks\relax
\def\floatpagepagefraction{1}
\def\textpagefraction{.001}

    \maketitle

    \abstract{
Sustainable forest management requires handling uncertainty introduced from various sources, considering different conflicting economic, environmental, and social objectives, and involving multiple decision-making periods. 
This study proposes an interactive and intuitive decision-support approach for sustainable, robust forest harvest scheduling in multiple periods in a short-term (6-12 months) planning horizon. 
The approach includes a novel multi-scenario multiobjective mixed-integer optimization problem that allows forest planners to separately study the trade-offs between demand satisfactions for multiple assortments in different planning periods. Moreover, it provides an intuitive robust analysis to support forest planners in dealing with uncertainty and investigating potential variations of the outcomes as the consequences of uncertainty in tactical forest planning problems.
We validate the proposed decision-support approach in a Swedish case study with 250 forest stands, three assortments (pine, spruce, deciduous trees), and a twelve-month harvest planning horizon. We demonstrate how the proposed approach supports a forest practitioner in trade-offs and robust analyses and finding the most preferred robust solution. 
    }
    
   \textbf{keywords:}
Forest harvest scheduling, 
Multiobjective optimization,
Multi-criteria decision-making,
Scenario planning,
Deep uncertainty,
Robustness analysis,
Sustainable forest planning

%    \textcolor{red}{Possible targeted journals: \\
%    - Journal of Environmental Management,
%    - Forest Policy and Economics,%\\
%    - International Journal of Forest Engineering, %\\
%    - Forest Ecology and Management,
    %- Canadian Journal of Forest Research,%\\
%    }

\section{Introduction}
We study problems of sustainable, robust forest harvest scheduling in multiple periods in a tactical (12-month) planning horizon with multiple objectives and uncertain timber volumes for various assortments. At this tactical harvest scheduling level, the objectives are to minimize the differences between timber production and demand for each assortment in every period. This is to reflect that industry losses occur for any mismatches (both shortfalls and surpluses) between timber deliveries and demand, thus the goal of forest planners is to be as precise as possible to fulfill the demand. While the industrial demand for each assortment in every period (month) is often assumed to be certainly known, the exact timber volume for each assortment is still uncertain, which can significantly affect the precision of the harvest schedule. Even gathering new data through conducting a forest inventory before harvest scheduling cannot wholly remove existing uncertainty in available timber volumes. Since fitting the proper statistical distribution might not always be possible, this classifies as \emph{deep uncertainty} \citep{Walker2013a, Shavazipour2019, Shavazipour2022forest}. 

In deep uncertainty, finding an optimal plan for an average or the most probable scenario (also called \emph{nominal} or \emph{base-line} scenario) may lead to failure \citep{Shavazipour2019}. Instead, identifying a robust plan that works relatively well in a broader range of scenarios is recommended, which requires specific uncertainty and robustness analyses \citep{Lempert2006, Shavazipour2019, Shavazipour2022forest}.
Although, in the last decades, various methodologies and tools have been developed for dealing with deep uncertainty in other disciplines \citep{Lempert2006, Kasprzyk2013, Herman2015, Quinn2017, Shavazipour2021a, shavazipour2023}, only a few studies consider handling deep uncertainty in the forest management context, mainly in strategic planning \citep{Yousefpour2016, Radke2017, Augustynczik2019, Radke2020, Horl2020, Shavazipour2022forest}. 
To the best of our knowledge, similar tools have yet to be developed or tested in tactical harvest scheduling problems. 

Besides the weakness of coping with deep uncertainty, traditional tools developed for harvest scheduling only consider aggregated demand satisfaction for different assortments across all harvesting periods (e.g., \citep{ronnqvist2003optimization, bredstrom2010annual}), which can hide industry losses caused by fluctuations between assortment productions and differences between periods. For example, when minimizing the deviation between the aggregating timber productions in different periods and the total annual demand, monthly demand satisfaction for a specific product is significantly ignored. Indeed, production shortfalls in multiple months for various products can be neutralized with significant extra production of a single product in a single month, which is not the aim of any planner. Therefore, the objective must be to minimize the production-demand deviation for all products in every period instead of minimizing the aggregated differences. This leads the planner to consider multiple objectives and perform a trade-off analysis of conflicting plans concerning demand satisfaction of various products and different periods. A multiobjective aspect of tactical forest harvest scheduling has also received little attention in the literature, mainly due to the extra complexity and computational expenses it will bring to the problem. 

Conflicting objective functions in an optimization problem lead to the situation of having numerous so-called Pareto optimal (or non-dominated) solutions instead of a single optimal one, which brings the need for a domain expert (e.g., a planner or decision-maker) to pick the final solution amongst various compromises. Since the set of all Pareto optimal solutions remains to be generated by solving mathematical optimization models, extra support is required for the planner to investigate the trade-offs and consequences of each Pareto optimal (or candidate) solution and identify the most suitable/preferred one for practical implementation. 
Over the years, various decision support based on multiobjective optimization methodology has been developed for different forest planning problems, e.g., \citet{Kangas2001, toth2009finding, Monkkonen2014, Kangas2015, Trivino2017, Eyvindson2018a, marques2021participatory, marques2021building}.

Adding deep uncertainty (and multiple scenarios) to a multiobjective optimization problem further complicates the planner's decision-making task \citep{Shavazipour2021b}. 
As mentioned earlier, different approaches like many-objective robust decision-making (MORDM) have been developed to provide support for trade-offs and uncertainty analyses for complex multi-scenario multiobjective optimization problems (e.g., \citet{Kasprzyk2012, Watson2017, eker2018, Shavazipour2019, Shavazipour2021a, shavazipour2024}). However, these approaches are not applied in forest harvest scheduling problems. In this study, we use one of the latest variants of MORDM: multi-scenario multiobjective robust optimization, called multi-scenario MORO, \citep{Shavazipour2021a, shavazipour2024} that proved its efficiency in generating more robust solutions by considering various scenarios already in the multiobjective optimization model.  

All uncertainty handling, solving multiobjective and multi-period planning problems, and modeling uncertainty for robust optimization and its practical usage in large-scale mixed-integer forestry planning problems are listed by \citet{ronnqvist2015operations} as open questions and research challenges in forest management. 

In this study, we propose a decision-support approach to address the literature gaps mentioned above. More specifically, the main contributions of this study can be summarized as follows:
\begin{itemize}
    \item To the best of our knowledge, this is the very first study that proposes a decision support approach for sustainable robust tactical harvest scheduling that
    allows forest planners to separately study the trade-offs between demand satisfactions for multiple assortments in different planning periods. 

    \item This is also the first study that proposed and tested a robust analysis to support forest planners in dealing with deep uncertainty and study the potential variations of the outcomes as the consequences of uncertainty in tactical forest planning problems. 

    \item We extend the usage of MORDM (particularly the multi-scenario MORO) in the forest management context.

    \item We introduce an efficient way of interacting with the decision-maker (forest planner) and solving a large-scale mixed integer multiobjective optimization problems with too many objective functions under deep uncertainty. 
\end{itemize}

In the next section, we introduce methodological elements including multi-scenario multiobjective optimization, multi-scenario MORO, and robustness analysis; and demonstrate the proposed interactive decision support approach for harvest scheduling. In Section 3, we frame the decision problem through a case study with descriptions of the study area and source data, the uncertainty and scenario generation, the mathematical problem formulation, and the solution method. In Section 4, we present and interpret the case study results describing how the proposed decision support can be applied in practice. In Section 5, we further discuss the feedback received from a practitioner in our case study experiment as well as other practical matters, such as computation expenses. Section 6 concludes the paper and raises future research directions.  

\section{Methodology}
\subsection{Multi-scenario multiobjective optimization problem}\label{msmoro}

A general form of a multi-scenario multiobjective optimization problem with $k \geq 2$ objective functions and $s \geq 2$ scenarios can be formulated as follows \citep{Shavazipour2021a}:
\begin{equation}
\begin{array}{rll}
\mbox{minimize}  & \mathbf{F}({\x})={\left(f_{1p}({\x}),  \dots, f_{kp}({\x}) \right)}, & p =1, \dots, s, \\ 
\mbox{subject to}  & h({\mathbf{x}}) = 0, \\
                   & g({\x}) \leq 0, & \\
\end{array}\label{eq:sbmop}
\end{equation}
where 
$f_{ip} $ denotes the $i$-th objective function under the conditions of scenario $p$ and $\x = (x_1, \dots, x_n)^T$ is a feasible solution represented by a vector of decision variables that
satisfy the equality ($h({\mathbf{x}})=0$) and inequality ($g({\mathbf{x}})=0$) constraints. The set of all feasible solutions forms the \emph{feasible region} in the \emph{decision space} $\Re^n$. 

The image of a feasible solution $\x$ conditional to a scenario $p$ in the \emph{objective space} $\Re^k$ is given by the \emph{objective vector} $\mathbf{z}_p =(f_{1p}({\x}),  \dots, f_{kp}({\x}))^T$. 
A \textit{Pareto optimal} solution in scenario $p$ is a feasible solution $\x^*$ which is strictly better than any other feasible solution $\x$ in at least one of the objective functions and at least one scenario and as good in the other objective functions and scenarios, i.e., satisfies 
$f_{ip}(\x^*) \leq f_{ip}(\x)$ for all $i,p$ and $f_{jq}(\x^*) < f_{jq}(\x)$ for at least one $j$ and one $q$ \citep{Shavazipour2021b}.
The image of the set of all Pareto optimal solutions in the objective space is also known as a \emph{Pareto front}. A Pareto optimal solution that satisfies the preferences of the decision-maker(s) regarding all  $k \times s$ objective-scenario combinations (also called meta-objectives) is called a \emph{preferred} solution. 

By solving all single-scenario single-objective optimization problems, the best possible value for each objective function in each scenario can be calculated. These values are the components of the so-called \emph{ideal vector} denoted by $\mathbf{z}^{\text{ideal}} = (z_{11}^{\text{ideal}}, \dots, z_{ks}^{\text{ideal}})^T$. Similarly, the \emph{nadir} vector $\mathbf{z}^{\text{nadir}} = (z_{11}^{\text{nadir}}, \dots, z_{ks}^{\text{nadir}})^T$ is composed by the worst possible value for each objective function in each scenario. In contrast to calculating the ideal values, computing the nadir values is not straightforward and usually requires approximation \citep{miettinen1999}.

\emph{Scalarization functions} have been widely used in the literature to incorporate the preferences of the decision-maker(s) and transform a multiobjective optimization problem into a single-objective problem. The optimal solution to the corresponding single-objective problem is a Pareto optimal solution to the original multiobjective problem \citep{miettinen1999, Miettinen2002, Ruiz2009classification}. The usage of scalarization functions (more specifically, the so-called achievement scalarizing function \citep{Wierzbicki1986}) has been recently extended to solve multi-scenario multiobjective optimization problems \citep{Shavazipour2020, Shavazipour2021a}. By using the multi-scenario version of an achievement scalarizing function (MS-ASF) suggested in \cite{Shavazipour2021a}, the problem \eqref{eq:sbmop} is scalarized to the following:
\begin{equation}\label{eq:sbasf}
\begin{array}{ll}
   \mbox{minimize}  & \left(\max\limits_{i=1,\dots,k,\, p=1,\dots,s} [w_{ip}(f_{ip}(\x) - \overline{z}_{ip})]\right) + \epsilon \sum\limits_{i=1}^k \sum\limits_{p=1}^s  w_{ip}(f_{ip}(\x) - \overline{z}_{ip})\\
  %&\\
     \mbox{subject to} &  h({\mathbf{x}}) = 0, \\
                   & g({\x}) \leq 0,  \\
\end{array}
\end{equation}
where the preference of the decision-maker(s) is given in terms of a desirable value (or \textit{aspiration level}) $\overline{z}_{ip}$, $z^{ideal}_{ip}\leq \overline{z}_{ip}\leq z^{nadir}_{ip}$, for each objective function in every scenario. The weight $w_{ip}$ describes the corresponding \emph{importance weight} that can be used, e.g., for normalization purposes. To guarantee Pareto optimality with respect to (\ref{eq:sbmop}) of the solution to (\ref{eq:sbasf}), a so-called augmentation term has been added, with $\epsilon$ a small positive scalar (we refer the reader to \cite{miettinen1999} for details of the role of an augmentation term in scalarization functions). Combining the aspiration levels for all objectives in all scenarios constructs a vector of desired values known as a \emph{reference point}. 
Different Pareto optimal solutions to (\ref{eq:sbmop}) can be generated by changing the reference points and/or the importance weights in (\ref{eq:sbasf}). % although sometimes the same Pareto optimal solution may be associated with multiple reference points or weights.

Scalarization functions can be applied \emph{a priori}, \emph{a posteriori}, or \emph{interactively} based on when the decision-maker provides their preferences \citep{miettinen1999}.
In this paper, an interactive method for solution generation might not be the best choice because of the problem's high dimensionality and the time-consuming calculations in the solution process. Therefore, we use the scalarization functions (\refeq{eq:sbasf}) in an a posteriori approach and solve the scalarized problem multiple times to pre-generate Pareto optimal solutions before the decision-making phase. Only the importance weights are modified between different solves; following \cite{Shavazipour2020}, we use the same aspiration levels because of the specific form of the problem and a large number of objective functions. Despite the a posteriori approach for the solution generation phase, the subsequent decision-making process, including the robustness and trade-off analyses, is interactive, similar to the multi-scenario MORO proposed in \cite{Shavazipour2021a}. 

\subsection{Multi-scenario multiobjective robust optimization approach}
\begin{figure}[!htbp]
    \centering
    \includegraphics[width=.7\textwidth]{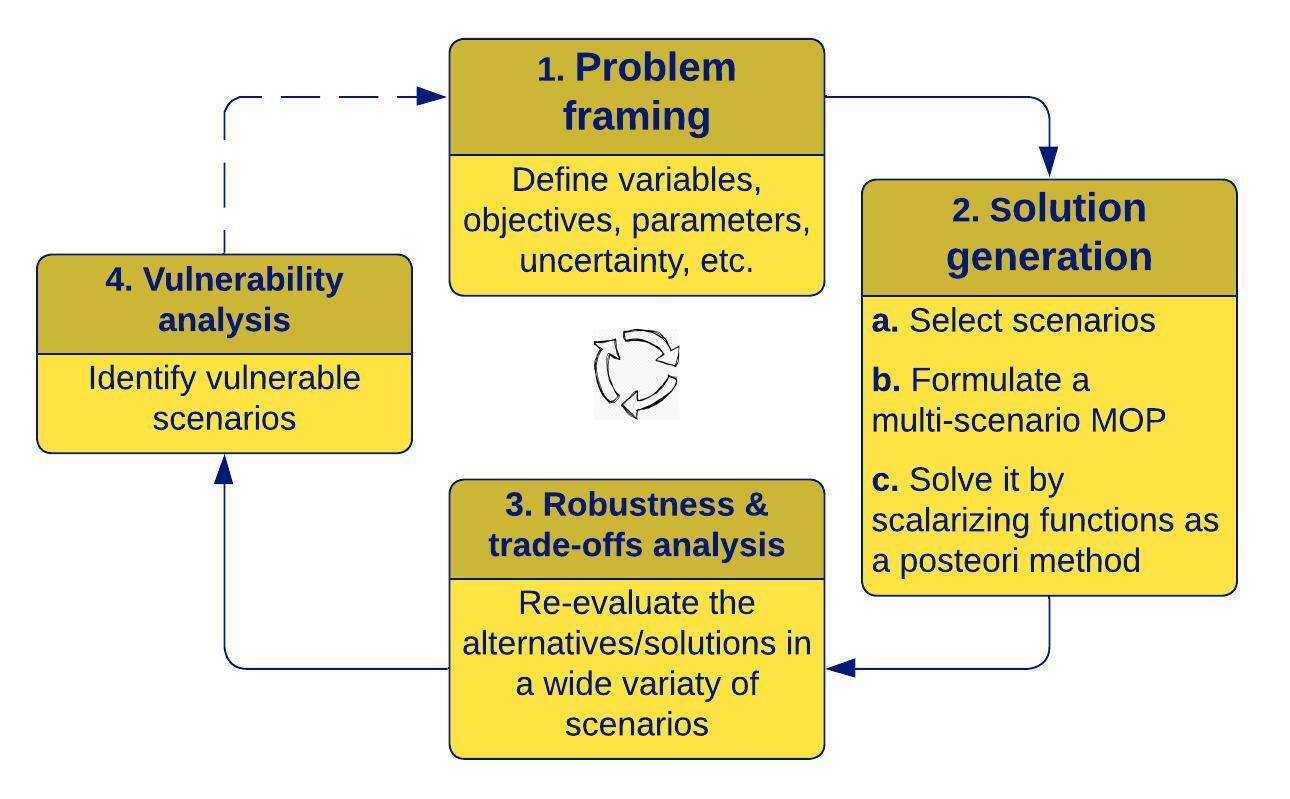}
    \caption{Schematic illustration of the multi-scenario multi-objective robust optimization approach (adapted from \cite{Shavazipour2021a}). }
    \label{fig:MSMORO}
\end{figure}

The multi-scenario MORO approach comprises four primary steps to iteratively identify the most preferred robust solution based on the decision-maker's preferences. Figure~\ref{fig:MSMORO} describes these four steps and how they relate to each other---i.e., (1) problem framing, (2) solution generation, (3) robustness and trade-off analyses (each candidate solution re-evaluated in a broad range of plausible scenarios), and (4) vulnerability analysis (identifying the scenarios that cause poor performance). This procedure persists until the decision-maker finds the most preferred solution. A more detailed description of the multi-scenario MORO approach, as per used in this study, is given in Algorithm 1 (Table~\ref{tab:algoritm}).

\subsection{Robustness and trade-off analysis}\label{robust_measures}
Several robustness measures have been developed in the literature to evaluate the robustness of alternative solutions in environmental decision-making, where deeply uncertain parameters cause significant performance variation (see, e.g., \cite{mcphail2018robustness} for detailed comparison and suitability of using various robustness measures). In this study, we use a satisficing metric called \emph{domain criterion} \citep{Starr1963} as it (1) is intuitive and easy to understand by different kinds of decision-makers with different levels of risk aversion (as it is based on user-selected minimum acceptable performance thresholds); (2) utilizes a broad set of plausible scenarios; and (3) can be adjusted by the decision-maker to meet their requirements and preferences in ---i.e., the decision-maker can interactively adjust the minimum acceptable
performance threshold, reflecting their risk aversion levels, and learning the feasibility, robustness, and consequences of different solutions (see \citet{shavazipour2024} for an interactive framework proposed for setting these thresholds).

The domain criterion measure calculates, for each solution, the number/percentage of scenarios in which the current solution can meet the acceptable threshold set by the decision-maker(s). A value of $1$ corresponds to the current solution meeting the given criterion in all scenarios, while a value of $0$ corresponds to the given threshold not being satisfied in any scenario. Any robustness value between these two ($[0,1]$) indicates the fraction of scenarios in which the satisfactory threshold has been met. The higher robustness value for a criterion, the more scenarios meeting the acceptable threshold for that criterion. In other words, a solution can be said being more robust than another solution if it has a higher robustness value.

\begin{table}[h!]
\vspace{10pt}\noindent
\begin{tabularx}{\textwidth}{l X}\toprule
    \multicolumn{2}{c}{Algorithm 1. Multi-scenario MORO decision support approach} \\\midrule
    \multicolumn{2}{l}{1. Problem framing:} \\
        1.1. & Define decision variables, objective functions, parameters, uncertainty, etc. \\[5pt] 
    \multicolumn{2}{l}{2. Solution generation:} \\
        2.1. & Formulate the multi-scenario multiobjective optimization model.\\
        2.2. & Generate a cohort of scenarios based on the deeply uncertain parameters.\\
        2.3. & Calculate the ideal values for every objective function in all generated scenarios.\\
        2.4. & Solve the scalarized problem, e.g., \eqref{eq:sbasf}, several times to generate a representative set of Pareto optimal solutions.\\[5pt]
    \multicolumn{2}{l}{3. Robustness and trade-off analysis (main decision-making process):} \\
        3.1. & Re-evaluate each Pareto optimal solution in all generated scenarios for stress testing. \\
        3.2. & Present the variation of objective function values across the generated cohort of scenarios to the decision-maker. \\
        3.3. & Ask the decision-maker to set their domain criteria (see \ref{robust_measures}) for each objective function (the minimum acceptable value) to assess the robustness of a solution (it can be done iteratively, e.g., the way proposed in \citep{shavazipour2024}. \\
        3.4. & Compute the robustness scores for all generated solutions and show the results to the decision-maker. If the number of objective functions is large, a subset including the most important robustness scores could be shown. \\
        3.5. & Save as many solutions as the decision-maker likes for further investigation. \\
        3.6. & Display all robustness scores for the chosen solutions to the decision-maker. \\
        3.7. & Show the decision variables of the chosen solutions to the decision-maker for further analysis if so desired. \\
        3.8. & If the decision-maker is not happy with a solution as the final one, do one of the following: \\
            & (a) Go to step 3.4 if the decision-maker wants to see all the generated solutions. \\
            & (b) Go to step 3.3 if the decision-maker wants to change the domain criteria. \\
            & (c) Go to step 2.4 if the decision-maker feels more Pareto optimal solutions need to be generated. \\
            & (d) Go to step 1 if the decision-maker feels the problem formulation needs to be updated. \\[5pt]
    \multicolumn{2}{l}{Stop: The final solution is the one identified in step 3.8.} \\
    \bottomrule
\end{tabularx}
\caption{Steps of the proposed algorithm for decision support.}\label{tab:algoritm}
\end{table} 

\section{Case study}\label{case-study}

\subsection{Study area and source data}
As already mentioned in the introduction of the paper, we study the problem of robust harvest scheduling in multiple time periods in a tactical (12-month) planning horizon with multiple objectives and uncertain timber volumes for multiple assortments. The scheduling task is to decide when (in which month $t = 1,\ldots,n_T$) each of the forest stands should be harvested so that the resulting aggregated flows of different timber assortments are Pareto optimal and accepted by the planner. 

\begin{figure}[!h]
    \centering
    \includegraphics[width=\textwidth]{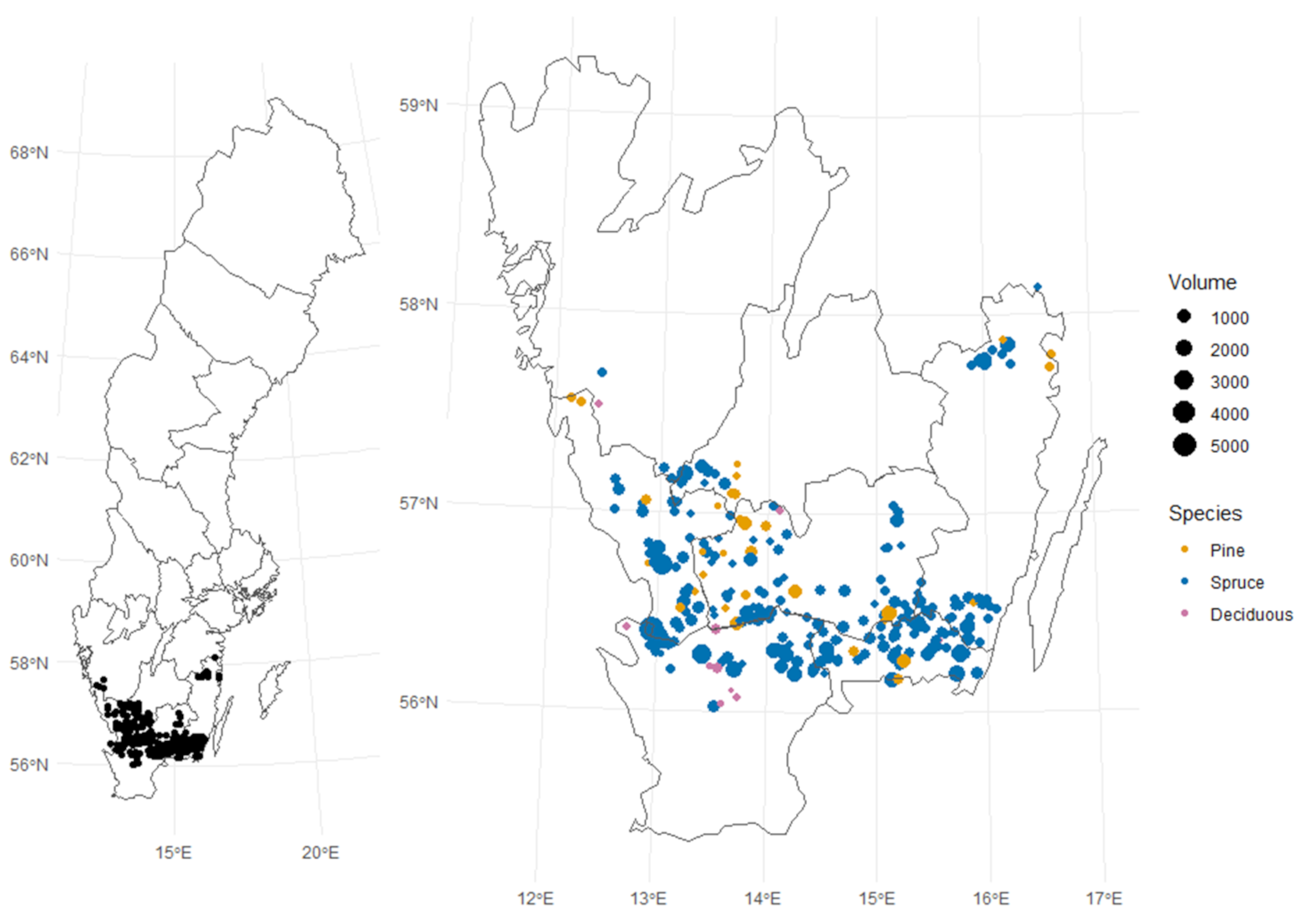}
    \caption{Forest stands in Southern Sweden included in the case study. The size of markers is proportional to the mean values of the volume of the dominating assortment.}\label{fig:study_area}
\end{figure}

The study area covers about $433$ hectares of forest in Southern Sweden (Figure~\ref{fig:study_area}). The forest is divided into $250$ forest management units (stands) which are dominated by spruce ($80~\%$ of stands), pine ($15~\%$), or deciduous trees ($5~\%$), however, all three tree assortments are to some extent available in all stands. The size of the stands varies from a minimum of $0.47$ hectares to a maximum of $7.66$ hectares, with an average of $1.73$ hectares. Figure~\ref{fig:study_area} illustrates how the stands are distributed within the landscape. The forest stand data required for the proposed decision support are the mean value and standard deviation of the volume of each tree assortment in every forest stand. Such pre-harvest data was simulated for the case study area using methods described in \cite{ene2021}. 

\subsection{Uncertainties and scenario generation}\label{uncertainty}
Three independent sources of uncertainty have been considered for each stand $j$, $j=1,\ldots,n_S=250$, resulting in a number of uncertain factors totaling to 750. The uncertainties are related to the available volume of the $n_A=3$ assortments $a=$~pine, spruce, deciduous. Since the exact volumes can be known only after harvesting, the harvest scheduling decisions are classified as ``here-and-now'' decisions---i.e., the decisions need to be made and implemented before realizing the actual volumes. 

The available information for the uncertain volume $V_{aj}$ of assortment $a$ in stand $j$ is a normal distribution with mean value $\mu_{aj}$ and standard deviation $\sigma_{aj}$. We use the interval $[\mu_{aj} - \sigma_{aj}, \mu_{aj} + \sigma_{aj}]$ for each assortment $a$ and stand $j$ to limit the scenario generation, and consider any combination of values within these 750 intervals as a plausible scenario. The scenario consisting of all mean values is used as the nominal scenario $p=2$. The scenario consisting of all lowest possible volumes (or 0 for those $(a,j)$ for which $\mu_{aj}-\sigma_{aj} < 0$) is referred to as the worst-case scenario and denoted $p=1$, and the scenario consisting of all highest possible volumes is referred to as the best-case scenario and denoted $p=3$ (whether these scenarios would actually be considered worst-case or best-case is not sure nor relevant). 

The three scenarios, $p=1,2,3$, are used when solving the multi-scenario MORO formulation of the harvest scheduling problem (defined in \eqref{eq:msp} below). In addition to those, a randomly generated, large cohort of scenarios is used for the robustness analysis. More specifically, we create 997 scenarios by making 997 realizations (drawn uniformly) for each of the 750 uncertain volumes $V_{aj}$, $\forall a,j$. This gives a total number of scenarios of 1000.  
The uncertain factors and scenarios are summarized in Table~\ref{tab:scenarios}. 

\begin{table}[h!]
\vspace{10pt}\noindent
\begin{tabularx}{\textwidth}{l l l l l}\toprule
    \multicolumn{5}{c}{Deeply uncertain factors and scenarios} \\\midrule
    Notation & Description & Range & Scenario name & Scenario value \\\midrule
    \multirow{3}{*}{$V_{aj}$} & Volume of & \multirow{3}{*}{$[\mu_{aj} - \sigma_{aj}, \mu_{aj} + \sigma_{aj}]$} & $p=1$ (worst-case) & $\max(\mu_{aj} - \sigma_{aj},0)$ \\ 
    & assortment $a$ & & $p=2$ (nominal) & $\mu_{aj}$ \\
    & in stand $j$& & $p=3$ (best-case) & $\mu_{aj} + \sigma_{aj}$ \\
    \bottomrule
\end{tabularx}
\caption{The deeply uncertain factors and the three scenarios considered when solving the multi-scenario MORO formulation of the harvest scheduling problem. Parameters $\mu_{aj}$ and $ \sigma_{aj}$ are the mean value and standard deviation for the volume of assortment $a$ in stand $j$.}\label{tab:scenarios}
\end{table}

\subsection{Mathematical problem formulation}
Following the general form of the multi-scenario multiobjective optimization problem \eqref{eq:sbmop}, the robust harvest scheduling problem is formulated as follows:
\begin{equation}\label{eq:msp}
\begin{aligned}
& \text{minimize}   && \mathmakebox[0pt][l]{\mathbf{F}(\mathbf{x}) = \left(f_1(\mathbf{x}), f_2(\mathbf{x}), \ldots, f_k(\mathbf{x})\right) } \\
& \text{subject to} && \sum_{t=1}^{n_T} x_{jt} \leq 1, &&  j=1,...,n_S \\
&  &&  x_{jt} \in \{0,1\},  &&  j=1,\ldots,n_S,\enskip t=1,\ldots,n_T
\end{aligned}
\end{equation}
where the $i$th objective function $f_i(\mathbf{x})$ is defined as 
\begin{equation*}
    f_i(\mathbf{x}) = \left| \left(\sum_{j=1}^{n_S} V_{ajp} x_{jt}\right) - D_{at} \right|,
\end{equation*}
with indices $(a,t,p)$ determined uniquely from $i$, $V_{ajp}$ denoting the volume of assortment $a$ in stand $j$ in scenario $p$, and $D_{at}$ denoting the assortment demand in period $t$. Expressed in words, the $k=n_A \times n_T \times s$ objective functions quantify the assortment-specific deviations of the total harvested volume to the demand under the conditions given by the scenarios. The decision variable $x_{jt}$ is defined as 
\begin{equation*}
    x_{jt} = \begin{cases} 
            1, & \text{if stand $j$ is harvested in time period $t$, and} \\
            0, & \text{otherwise}.
             \end{cases}
\end{equation*} 
The inequality constraints in \eqref{eq:msp} thus ensure that each stand is scheduled for harvesting at most once. In our case study, the number of objective functions is $k = 3 \times 12 \times 3 = 108$, and the number of decision variables is $250 \times 12 = 3000$.

\subsection{Solution method}\label{solution}
To solve the high-dimensional nonlinear multiobjective optimization problem \eqref{eq:msp}, we first linearize the problem by introducing the artificial variable $\phi_{atp}$, resulting in the following formulation:

\begin{equation}\label{eq:mslmop}
\begin{aligned}
& \text{minimize} && \mathmakebox[0pt][l]{\phi_{atp}, \quad a=1,\ldots,n_A, \enskip t=1,\ldots,n_T, \enskip p=1,\ldots,s} \\
& \text{subject to} && \phi_{atp} \geq \left(\sum_{j=1}^{n_S} V_{ajp} x_{jt}\right) - D_{at}, && \forall a,t,p \\
&&& \phi_{atp} \geq D_{at} - \left(\sum_{j=1}^{n_S} V_{ajp} x_{jt}\right), && \forall a,t,p \\
&&& \sum_{t=1}^{n_T} x_{jt} \leq 1, &&  \forall j \\
&&& x_{jt} \in \{0,1\}, &&  \forall j,t \\
&&& \phi_{atp} \geq 0, && \forall a,t,p
\end{aligned}
\end{equation}

The multiobjective optimization problem \eqref{eq:mslmop} is then transformed into a single-objective problem through the multi-scenario achievement scalarizing function as defined in \eqref{eq:sbasf}, resulting in the following formulation: 

\begin{equation}\label{eq:msasf}
\begin{aligned}
& \text{minimize}  && \mathmakebox[0pt][l]{\left(\max\limits_{a,t,p} \left[w_{atp}(\phi_{atp} - \overline{z}_{atp})\right]\right) + \epsilon \sum_{a=1}^{n_A} \sum_{t=1}^{n_T} \sum_{p=1}^{s} w_{atp}(\phi_{atp} - \overline{z}_{atp})} \\
& \text{subject to} && \phi_{atp} \geq \left(\sum_{j=1}^{n_S} V_{ajp} x_{jt}\right) - D_{at}, && \forall a,t,p \\
&&& \phi_{atp} \geq D_{at} - \left(\sum_{j=1}^{n_S} V_{ajp} x_{jt}\right), && \forall a,t,p \\
&&& \sum_{t=1}^{n_T} x_{jt} \leq 1, &&  \forall j \\
&&& x_{jt} \in \{0,1\}, &&  \forall j,t \\
&&& \phi_{atp} \geq 0, &&   \forall a,t,p
\end{aligned}
\end{equation}

As a final step, the nonlinear single-objective optimization problem \eqref{eq:msasf} is again linearized by introducing the artificial scalar variable $\varphi$, resulting in the following mixed-integer linear single-objective optimization problem that can be solved by any mixed-integer linear solver available in the literature:

\begin{equation}\label{eq:lms-asf}
\begin{aligned}
& \text{minimize} && \mathmakebox[0pt][l]{\varphi + \epsilon \sum_{a=1}^{n_A} \sum_{t=1}^{n_T} \sum_{p=1}^{s} w_{atp}(\phi_{atp} - \overline{z}_{atp})} \\
& \text{subject to} && \varphi \geq w_{atp}(\phi_{atp} - \overline{z}_{atp}), && \forall a,t,p \\
&&& \phi_{atp} \geq \left(\sum_{j=1}^{n_S} V_{ajp} x_{jt}\right) - D_{at}, && \forall a,t,p \\
&&& \phi_{atp} \geq D_{at} - \left(\sum_{j=1}^{n_S} V_{ajp} x_{jt}\right), && \forall a,t,p \\
&&& \sum_{t=1}^{n_T} x_{jt} \leq 1, && \forall j \\
&&& x_{jt} \in \{0,1\}, && \forall j,t \\
&&& \phi_{atp}, \varphi \geq 0, && \forall a,t,p
\end{aligned}
\end{equation}

\section{Results}\label{results}
%In this section, we adapt the steps of the multi-scenario MORO approach, shown in Fig.~\ref{fig:MSMORO}, to be used as a decision support for our forest harvest scheduling problem. 
%In this section, 
%We also provide empirical evidence about the usefulness of the proposed decision support via experiment with a decision-maker.

\subsection{Applying the multi-scenario MORO approach}\label{sec:apply_MSMORO}
We now adapt and apply the steps of the multi-scenario MORO approach detailed in Algorithm 1 (\ref{tab:algoritm}) to the harvest scheduling problem. The initial problem framing (step 1) and problem formulation (step 2.1) steps have already been described in Section~\ref{case-study}. 
%\url{Step_1.} \url{Problem_framing}\label{step1}
%As described in Section~\ref{case-study}, the decisions to be made are the harvest schedule for each stand (which stands to be harvested at each period). The objectives are to minimize any deviations from the demand for each of the three assortments (pine, spruce, and deciduous) during each period (month). The sources of uncertainty are the available assortment volumes ($m^3$) in $250$ stands.

\subsubsection*{Step 2. Solution generation}\label{step2}
%\url{Step_2.} \url{Solution_generation}
We generate solutions to the multi-scenario multiobjective harvest scheduling problem (i.e., harvest schedules) by solving the optimization problem \eqref{eq:lms-asf}, which is the linearized and scalarized version of the original formulation \eqref{eq:msp}. With the $s = 3$ scenarios, $p=1,2,3$ (see Section~\ref{uncertainty}), $n_T = 12$ periods, and $n_A = 3$ assortments, the original formulation becomes a 108-objective mixed integer optimization problem. With $n_S = 250$ forest stands, the final linearized and scalarized version \eqref{eq:lms-asf} becomes a mixed integer problem with $3000$ binary and $109$ %216$ 
real decision variables, and $476$ constraints. The problem is solved with Gurobi optimizer 10.0.2. The code and data are freely accessible at \url{-Link_to_be_added_after_the_acceptance (we_will_also_provide_it_to_the_reviewers_upon_request)}. 

\begin{figure}[h!]
    \centering
    \includegraphics[width=1.01\textwidth]{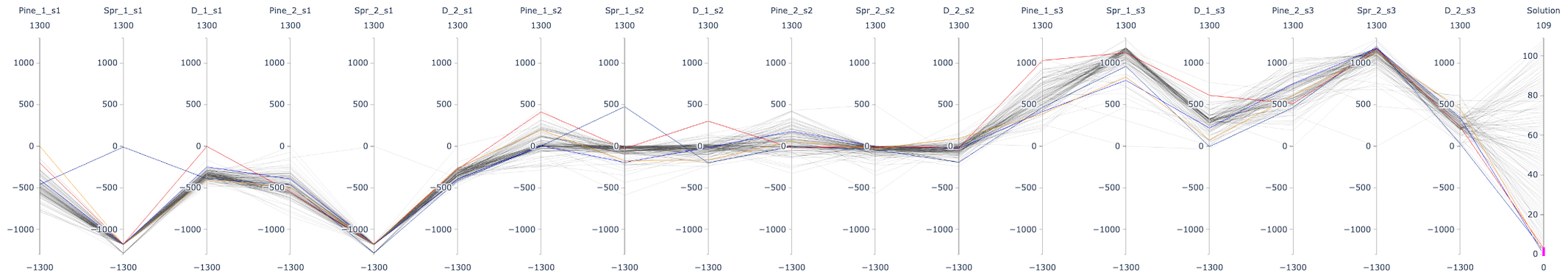}
    \caption{Parallel plot of objective values for the 109 generated solutions to Problem~\eqref{eq:lms-asf} (i.e., Pareto optimal solutions to the original formulation in Problem~\eqref{eq:msp}). For readability, only the objectives for the first two time periods are shown, resulting in showing 18 of the 108 objectives ($a = $ pine, spr(uce), d(eciduous), $t = 1,2$, and $p = 1,2,3$). As examples, the three solutions that provide the lowest deviation (objective value closest to zero) for respectively pine (yellow), spruce (blue), and deciduous (red) in time period $t=1$ of scenario $p=1$ have been highlighted. The right-most axis keeps track of the solution number.} \label{fig:109sol-3per}
\end{figure}

To generate different Pareto optimal solutions, Problem~\eqref{eq:lms-asf} needs to be solved multiple times. As mentioned in Section~\ref{msmoro}, we keep the aspiration levels $\overline{z}_{atp}, \forall a,t,p$, at the same values and only vary the importance weights $w_{atp}, \forall a,t,p$ of the objectives.
We generated in total 109 Pareto optimal solutions, consisting of the 108 solutions generated by giving an importance weight of 100 to each of the objective functions at a time (and the others a weight of 1), and the neutral solution generated by giving equal important weights `1' to all objective functions. Figure~\ref{fig:109sol-3per} shows all the 109 generated solutions in a parallel coordinate plot. Note that because of the high number of periods (thus objective functions), for the sake of simplicity and better readability, we only present the objective values for the first two time periods. The objective values for all twelve periods are available for further investigation in a later step in the decision-making process. Each poly-line in the plot describes the objective values for one solution. The intervals spanned by the poly-lines on each axis indicate that there are significant trade-offs between the objective functions. 

\subsubsection*{Step 3. Robustness and trade-off analysis}\label{step3}
\begin{figure}[h!]
    \centering
    \includegraphics[width=1.01\textwidth]{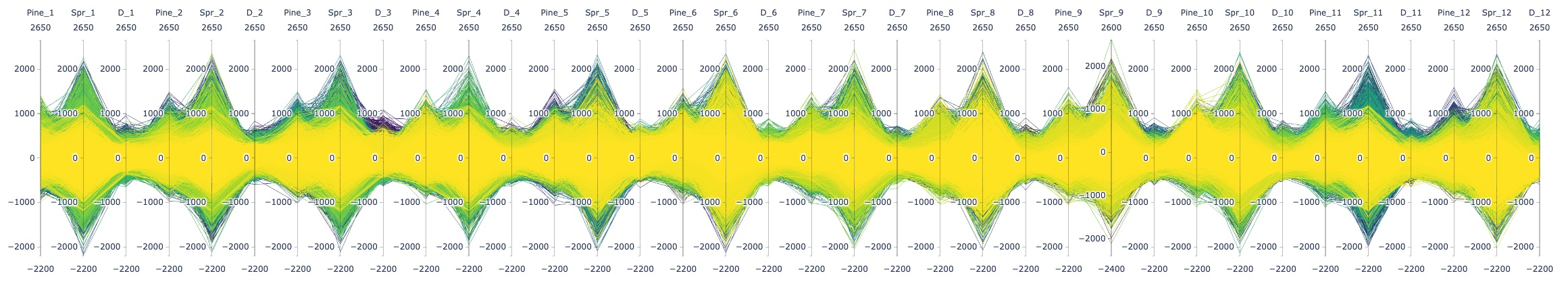}
    \caption{Parallel plot of the 36 objective functions evaluated for the 109 candidate solutions across the set of 1000 scenarios ($a = $ pine, spr(uce), d(eciduous), and $t = 1,\ldots,12$).}
    \label{fig:all108}
\end{figure}

%\url{Step _3.} \url{Robustness_and_trade-off_analysis} 
%As described in Section~\ref{robust_measures}, the domain criterion robustness measure is used in this study. 
To calculate the robustness of a candidate solution generated in \hyperref[step2]{Step 2}, the decision-maker first sets a minimum acceptable performance threshold for each of the 36 objective functions (36 is the number of objective functions if the scenario index is removed). The thresholds are part of the domain criterion measure (see Section~\ref{robust_measures}). Then, the candidate solution is re-evaluated across the set of $1000$ generated scenarios (see Section~\ref{uncertainty}) to check the number/percentage of scenarios in which the thresholds are met. 
Figure~\ref{fig:all108} shows the variations in objective function value obtained for the 109 candidate solutions when evaluated in the 1000 scenarios. Different colors distinguish solutions. %Robustness analysis will be discussed later when we describe our interactive decision-making experiment with a decision-maker in Section \ref{interactive}.  

\subsubsection*{Step 4. Vulnerability analysis}\label{step4}
%\url{Step_4.} \url{Vulnerability_analysis} \\
Following \citet{Shavazipour2021a}, the vulnerability/feasibility check (where the solutions cannot satisfy the domain criterion in an objective function across various scenarios) has been done by calculating the ideal points for all $36$ objective functions in all $1000$ scenarios (i.e., the best possible value for each objective function in each scenario). To compute the ideal points in our case study, one needs to solve $36\times1000 = 36$ $000$ single-scenario single-objective optimization problems. Nonetheless, the solution process for a single problem is quite fast (less than a second), and calculations only need to be done once. Furthermore, since, in our case study, the demand for each assortment is equal in different periods, we only need to solve the problem in one period (i.e., in total, $3000$ optimization problems need to be solved instead of $36$ $000$), as the optimization problems would be exactly the same in the other periods.
%Using a MacBook Pro with an Apple M2 Pro CPU and 16 GB RAM, all the calculations (solving $3000$ optimization problems) took about 10 minutes with the Gurobi optimizer 10.0.2.

The minimum and maximum ideal values for each objective function (deviation from the wood production and demand in a period) across all $1000$ scenarios have been described in Table \ref{tab:ideals}---i.e., the best possible values for each objective function in the best- and worst-case scenario. 
Studying the variation of the ideal values for each assortment across the scenarios shows that it is possible to relatively fulfill the demand in any scenario with a very small deviation from the demand if only one assortment and one period are considered. For example, the demand deviation for spruce in at least one scenario is $0.049$, meaning that even the optimal solution (specifically found for spruce in a period) is not exactly zero under conditions of at least one scenario.

\begin{table}[!h]
%\resizebox{\textwidth}{!}{
%\thispagestyle{empty}
\centering%
\begin{tabular}{cccc}
  \toprule
\multirow{2}{*}{} &  \multicolumn{3}{c}{Objective functions}  \\  \cline{2-4}
 & Pine & Spruce &  Deciduous \\ \hline
 
min &	4.229e-08 &	7.878e-07 &	4.143e-09	\\ \hline
max &	 0.013    &	0.049     &	 0.0043		\\ 
\bottomrule
\end{tabular}
%}
\caption{Maximum and minimum values for each objective functions among the components of the ideal points across $1000$ randomly generated scenarios.}\label{tab:ideals}
\end{table}

\subsection{An interactive decision-making experience}\label{interactive}
We used a \emph{Jupyter Notebook} to implement the decision support prototype and interact with the decision-maker. The decision-maker in our experiment was an expert in forest management. 
We solve the mixed integer linear optimization problem \eqref{eq:lms-asf} with Gurobi optimizer 10.0.2.
The code and data are freely accessible at \url{-Link_to_be_added_after_the_acceptance (we_will_also_provide_it_to_the_reviewers_upon_request)}. 

\subsubsection*{Initial calculations}
%We follow the proposed Algorithm 1 and start the decision-making process 
In the initialization (pre-experiment) phase, after framing the problem, generating $1000$ random scenarios, calculating ideal values for each objective function in all scenarios, and $109$ Pareto optimal solutions, we re-evaluate all the solutions in all scenarios to investigate variations in objective functions across generated scenarios (possible ranges). Then we have all the required data to conduct the decision-making process. 

\subsubsection*{Decision-making process}
\begin{figure}[!h]
    \centering
    \includegraphics[width=1.01\textwidth]{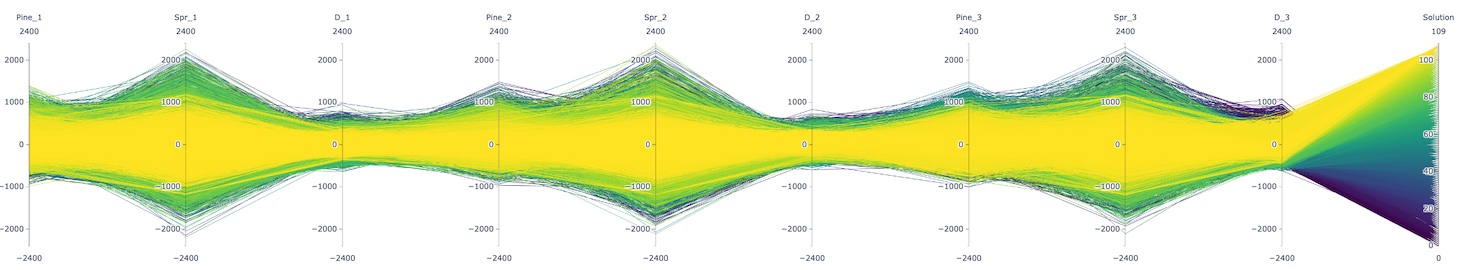}
    \caption{The possible ranges of objective functions across $1000$ randomly generated scenarios in the first three periods. Pine\_t: deviation from the demand for pine in period $t$, Spr\_t: deviation from the demand for spruce in period $t$, D\_t: deviation from the demand for deciduous in period $t$, $t=1,\dots,3$.
    }
    \label{fig:obj_variations}
\end{figure}

We started the decision-making process by giving the decision-maker an introduction to the problem, the optimization model, and the whole decision-making process. 
Then, the possible ranges of objective functions across $1000$ randomly generated scenarios were shown to the decision-maker to give an overall view of the variations in the outcomes caused by uncertainties (Figure~\ref{fig:obj_variations}). Note that to reduce the cognitive load and complexity, we only showed the results for the first three periods at this time based on the decision-maker requests. Another reason was that scheduling the first three periods was more important to be more robust, but the later schedules can be adapted later. In general, this process can be done for each quarter separately, which can also reduce uncertainty, leading to more adaptive decision-making. 

\subsubsection*{Iteration 1.}
In the next step, we represented the demand for each assortment in a period (demand for pine, spruce, and deciduous in each period are $1694$, $9227$, and $583$, respectively.) and asked the decision-maker to provide their preferences on the robustness measurement of the candidate solutions (domain criterion robust analysis). 
The decision-maker provided their preferences as the \emph{maximum acceptable percentage of deviation from the demand} for each assortment in each period. 
Therefore the \emph{domain criteria} were: 
\begin{itemize}
    \item Pine demand deviation $<$ $30\%$ 
    \item Spruce demand deviation $<$ $20\%$
    \item Deciduous demand deviation $<$ $30\%$
\end{itemize}
%$30\%$, $20\%$, and $30\%$ deviation from the demand for pine, spruce, and deciduous in each period, respectively. 
these preferences were the same in all periods. 

\begin{figure}[!h]
    \centering
    \includegraphics[width=1.01\textwidth]{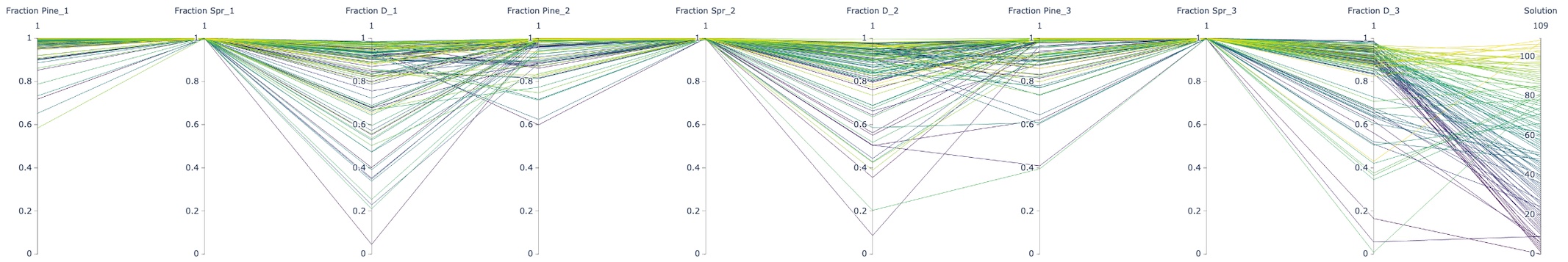}
    \includegraphics[width=1.01\textwidth]{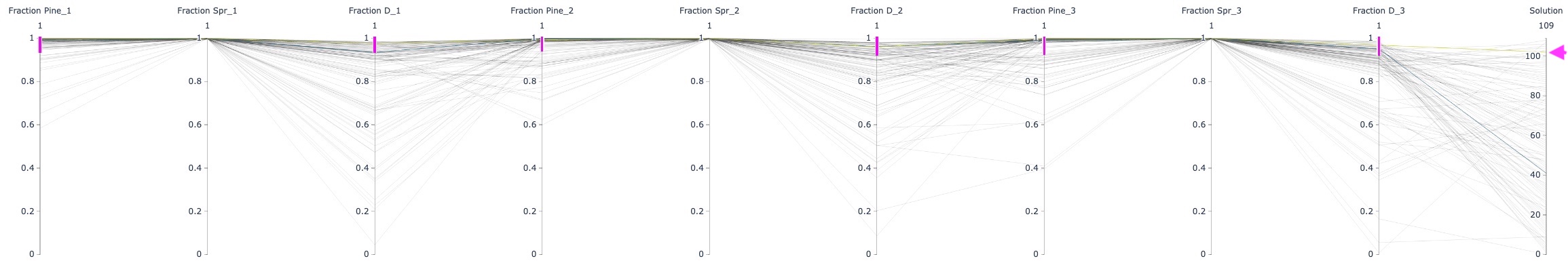}
    \caption{The robustness comparison of all $109$ generated solutions in different objective functions in the first three periods (top figure)), and the two most robust solutions in the first three periods (bottom).
    }
    \label{fig:robustness}
\end{figure}

Then, for each solution, we calculate the number of scenarios in which corresponding objective values satisfy each criterion mentioned above across all $1000$ scenarios. After further ranking and sorting the solutions based on their robustness in each criterion, we represented the robustness scores to the decision-maker as shown in Figure~\ref{fig:robustness}-top. Each line in this figure represents the robustness of a solution in different objective functions (demand deviations from the demand for pine, spruce, and deciduous in the first three periods). 
As requested by the decision-maker, we filtered out the solutions with less than $95\%$ robustness in the objective functions in the first three periods, as the bottom figure in Figure~\ref{fig:robustness} describes. As seen in this figure, only one solution (solution number `102') matched this constraint, followed by another one (solution number `42') with the robustness of $93.5\%$ for deciduous demand deviation in period 1, and higher than $95\%$ robustness in the others (within the first three periods). The decision-maker chose solution number $102$ for further investigation.  

\begin{figure}[!h]
    \centering
    \includegraphics[width=1.01\textwidth]{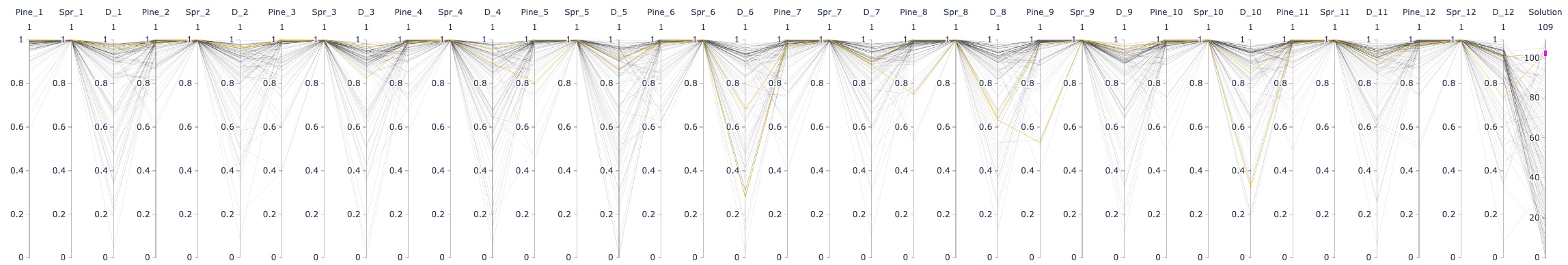}
    \caption{Robustness scores for solution number $102$ in all twelve periods.
    }
    \label{fig:robustness_102_12}
\end{figure}

Checking its robustness in other periods (periods 4-12) shows low robustness in some objective functions in the later periods, particularly for deciduous demand deviations in the sixth ($29\%$), tenth ($32\%$), and eighth ($64\%$) periods, see Figure~\ref{fig:robustness_102_12}.   
Therefore, the decision-maker wished to investigate to see if there were other options with better robustness scores than this solution (solution number `$102$') in the later periods but, perhaps, with less robustness in the first three periods.

\subsubsection*{Iteration 2.}
\begin{figure}[!h]
    \centering
    \includegraphics[width=1.01\textwidth]{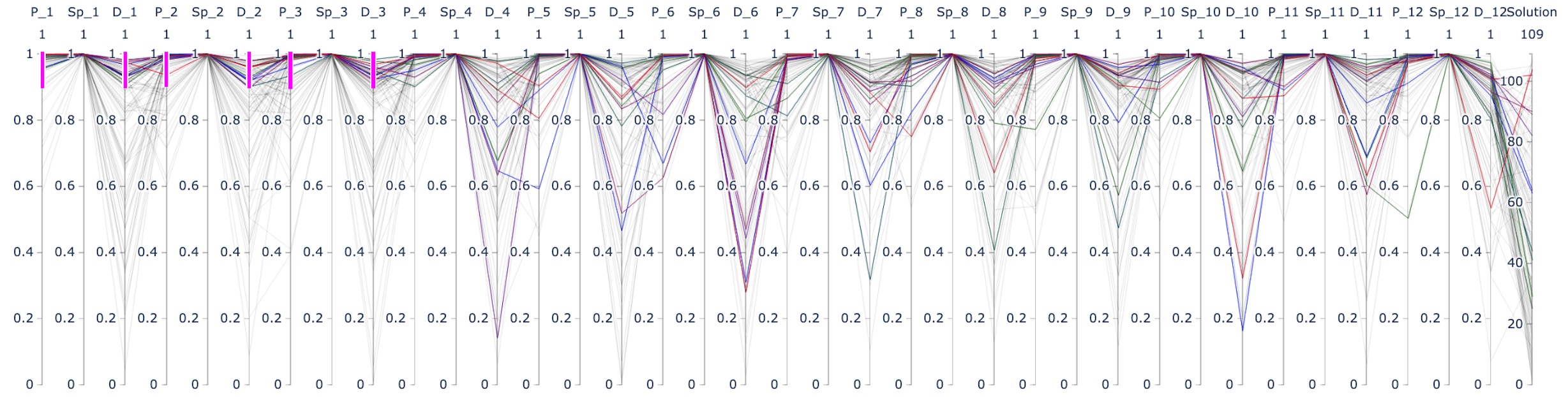}
    \includegraphics[width=1.01\textwidth]{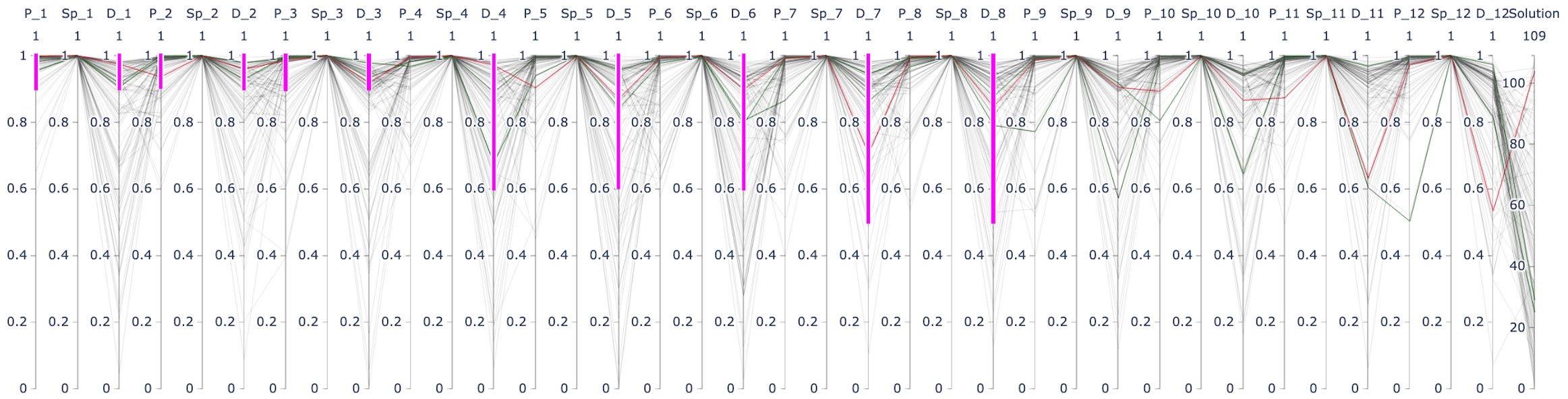}
    \caption{Solutions with the robustness of higher than $90\%$ in the first three periods (eleven solutions shown in the top figure), and the three candidate solutions with the robustness scores of higher than $60\%$ in the second quarter and $50\%$ in the last half of the year. In both figures, the robustness scores are presented in all twelve periods.
    }
    \label{fig:robustness_11-3}
\end{figure}

Therefore, in the second iteration, we widened the robustness restriction in the first three periods to $90\%$ and presented the eleven solutions found to the decision-maker. As seen in Figure~\ref{fig:robustness_11-3}--top, this time, we showed the robustness scores in all twelve periods. We then further eliminated the solutions with robustness scores of less than $60\%$ in the second quarter and less than $50\%$ in the last half of the year and displayed the remaining three solutions (number `$25$', `$29$', and `$104$') to the decision-maker, Figure~\ref{fig:robustness_11-3}--bottom. Solution number `$104$' is the most robust solution in the first half of the year, with robustness scores of higher than $87\%$. However, its robustness for deciduous demand deviation is the worst among these three ($53\%$). 
Solution number `$25$', $100\%$ meets the domain criteria for spruce demand deviations in every single period, but the ninth one ($99.6\%$). However, its robustness scores for other assortments are not competitive with the other two solutions, particularly in the last three periods.  
Finally, the robustness of solution number `$29$' has never gotten worse than $65\%$ in any period. 

Then, before making the final decision, for further investigation of the suitability of the candidate solutions, the decision-makers wanted to compare the corresponding decision variables (i.e., the stands to be harvested in each period) for these three solutions. Therefore, this information is shown to them as portrayed in Table \ref{tab:decisions}.

After analyzing the decision variables, the decision-maker realized that the number of stands to be harvested in some periods (namely period $11$ in solution `$104$' and period $10$ in solution `$29$') is too low compared to the other periods, meaning that there are probably larger stands, which may lead to some practical difficulty, particularly in cold and rainy periods, like period 10 and 11, when the uncertainty related to available working days and roads are high. 
Therefore, the decision-maker liked solution `$25$' the most, which has a more even distribution of the stands to be harvested with respect to their area. 
After the final checkout of the robustness of this solution, the decision-maker chose solution number `$25$' as the final solution.
As mentioned, this solution satisfies almost one hundred percent of the domain criterion for spruce demand deviations (spruce is the assortment that dominates the others) in every single period. 
Also, the robustness of this solution in the earlier periods was pretty high, and the lower robustness in the last quarter for deciduous (which has the lowest demand) and pine could be addressed by adjusting the plan later or sending the extra products left from earlier periods. 
Thus, the solution seemed fine with the decision-maker.

\begin{table}
    \centering
    %\resizebox{\columnwidth}{!}{%
    \begin{tabular}{|c|c|c|c|} \hline 
 & \multicolumn{3}{|c|}{Number of stands to be harvested}\\ \hline  
 Period& Solution 25&Solution 29 &Solution 104\\\hline \hline 
         1&  24& 20&26\\ \hline 
         2&  19& 21&24\\ \hline 
         3&   23& 26&26\\ \hline 
         4&  27&13&20\\ \hline 
         5&  25&20&19\\ \hline 
         6&   22&19&15\\ \hline 
         7&   20&23&20\\ \hline 
         8&   13&25&20\\ \hline 
         9&   15&21&17\\ \hline 
         10&   21&\textbf{\textcolor{red}{11}}&23\\ \hline 
         11&   17&24&\textbf{\textcolor{red}{8}}\\ \hline 
         12&   22&24&24\\ \hline
    \end{tabular}%
    %}
    \caption{Number of stands to be harvested in solutions number 25, 29, and 104, in period $t$, $t = 1, \dots, 12 $.}
    \label{tab:decisions}
\end{table}

\section{Discussion} 

\subsection{Feedback from the decision-maker}
As mentioned earlier, our decision-maker was a domain expert with several years of experience in tactical forest harvest scheduling. 
Here, we summarized the feedback we got from the decision-maker regarding their experiment with the proposed decision support prototype. 

\textbf{Uncertainty and robustness:} The decision-maker liked the overall view of the variations in the objective functions caused by uncertainty that was shown through the visualization of the possible ranges of objective functions across $1000$ randomly generated scenarios. They find it informative and useful to see how uncertainty can affect the outcome of a solution. 
It also helped them get a better idea about the complexity and their expectations, as well as not setting too optimistic domain criteria for robustness. 

\textbf{Robustness in the earlier periods:}
They also emphasize the importance of robustness in the earlier periods. The harvest schedule for the later periods can be updated later, after getting more information about the actual harvest volumes in the earlier periods. Usually, there is a bit of flexibility in late delivery of some parts and/or filling the shortcomings of harvest in one period with some extra harvested wood from the previous period. 

\textbf{Evenly distributed number of stands across the year:}
Moreover, they pointed out the importance of having evenly distributed stands scheduled across different periods and avoiding scheduling a few stands for harvesting in a single period, particularly in areas with limited accessibility (see, e.g., \cite{jonsson2017potential, westlund2024}). Indeed, to have a more robust harvest schedule, one also needs to consider the accessibility of the roads, particularly in rainy seasons, and consider more options for harvesting to avoid pauses in harvesting due to inaccessible stands. 
In similar cases, harvest planners prefer not to have a few stands in the harvest schedule in a single period to leave room for adjustment if required, e.g., the harvest team can harvest some other stands if they cannot access a specific stand, scheduled to be harvested, in some days/weeks. 
In fact, although it will bring more uncertainty and complexity to the problem, route accessibility in different periods should be considered in the optimization model \citep{westlund2024}, leading to more robust and adaptive planning that would be an interesting and important future research direction. 

\textbf{No such support exists despite the need:}
Finally, the decision-maker points out that, despite the need, a decision support tool, like the proposed one, is lacking in practice and encourages us to continue working on this line of research. 

\subsection{Computation cost}
All the calculations, in this study, have been done using a MacBook Pro with an Apple M2 Pro CPU and 16 GB RAM. All the optimization models are solved by Gurobi optimizer 10.0.2.
The main resources of the computation cost of the proposed decision support (and multi-scenario MORO in general) are (1) the solution generation step, % (in \url{Step_2}),
(2) stress test or solution re-evaluation, % (in \url{Step_3}), %3998.635487794876 seconds 
and (3) calculating the ideal values for each objective function in each scenario. % (in \url{Step_4}). 
All other calculations took no more than a few seconds and do not need to be considered.
Table \ref{tab:comp_cost} represents these sources, how much they cost in our case study, and the most effective parameters on the cost of each source. 

As described in this table, in total, all the calculations in our case study with $36$ objective functions, three scenarios in the mixed integer linear optimization model with $476$ constraints, $3000$ binary, and $109$ %216$ 
real decision variables, $1000$ scenarios for stress testing, and $109$ generated solutions took a bit more than three hours. 
Note that all of these calculations only need to be done once and before the decision process. So, they would not introduce any waiting time to the decision-maker in the decision-making process. 
However, for larger-scale problems (e.g., problems with a -higher number of stands, assortments, periods, and management options) and cases with a higher number of objective functions, scenarios, and solutions to be generated, this pre-decision-making process would possibly take longer time. 
The most effective parameter that causes more complexity and significantly increases the solution time in the multi-scenario multi-objective optimization models is the number of scenarios (and also the number of objectives) considered in these models. We refer the interested reader to \citet{Shavazipour2021a} for a detailed analysis of considering higher numbers of scenarios within the multi-scenario MORO approach and their computation costs. 

Also, calculating the ideal values for all objective functions in all scenarios and stress testing could be extremely time-consuming, particularly if tens of thousands or millions of scenarios are generated for stress testing of many solutions. 
As mentioned in Section \ref{step4}, because of the similar demand for each assortment in every period, we could twelve times reduce the number of optimization problems to be solved from $36$ $000$ to $3000$. Otherwise, it could take much more time than the current one. 
Nonetheless, as said, one is required to do all of these calculations once and before starting the decision-making process. Furthermore, all the calculations for ideal values and stress testing can be done in parallel, which can hugely reduce the total computation time if so desired. 

\begin{table}[!h]
    \centering
    \resizebox{\columnwidth}{!}{
    \begin{tabular}{c||c|c|l} 
             & Single evaluation& In total& Effective parameters\\ \hline \hline 
        \multirow{3}{*}{Solution generation} & \multirow{3}{*}{$\approx 60''$ per solution} & \multirow{3}{*}{$109'$} & number of objective functions,\\
        &  & &number of scenarios in the optimization model, \\ 
        &  & &number of solutions to be generated. \\ \hline 
        \multirow{2}{*}{Ideal values} & \multirow{2}{*}{$<1''$ per optimization} & \multirow{2}{*}{$10'$} & number of objective functions,\\
        &  & &number of scenarios generated for the stress testing, \\ 
        &  & &number of solutions to be generated. \\ \hline 
        \multirow{3}{*}{Stress testing} & \multirow{3}{*}{$\approx 37''$ per solution} & \multirow{3}{*}{$\approx 67'$} & number of objective functions,\\
        &  & &number of scenarios generated for the stress testing, \\ 
        &  & &number of solutions to be generated. \\ \hline
        % &  & &\\ \hline
    \end{tabular}
    }
    \caption{Main sources of the computation cost, their expenses in our case study, and the most effective parameters on the expenses.}
    \label{tab:comp_cost}
\end{table}

%\subsection{Limitations} mention in the conclusions

\section{Conclusions}
This study proposed a novel interactive decision-support approach for sustainable forest harvest scheduling under deep uncertainty using multi-scenario MORO \citep{Shavazipour2021a}. 
The sources of deep uncertainty considered here are the available volume ($m^3$) for three assortments in each stand: pine, spruce, and deciduous.
The objectives are to minimize any deviation from the demand for each assortment in every period. 
In the light of three scenarios (an average, the lowest, and the highest possible assortment volumes in each stand), three assortments, and twelve periods, we faced a mixed integer multiobjective optimization problem with $108$ objective functions, $476$ constraints, and $3109$ decision variables, that solved by a multi-scenario achievement scalarizing function several times to generate $109$ Pareto optimal solutions. 
All the generated solutions are re-evaluated (stress-tested) in an ensemble of $1000$ randomly generated scenarios to assess their robustness via domain criterion robustness measure. 
To the best of our knowledge, this is the first study proposing decision support for handling deep uncertainty and robust analysis in tactical %(short-term) 
forest harvest scheduling. A prototype of the proposed decision support has also been successfully validated with a domain expert, increasing the usefulness and potential usage in practice. 

In contrast to stochastic optimization models, which search for the optimal solution in a single average-case scenario, the proposed approach helps investigate a more comprehensive range of characteristics of forest harvest scheduling, trade-offs analysis, and variations in uncertain parameters. The proposed approach also supports decision-makers in stress-testing candidate solutions in multiple periods and various plausible scenarios and finding the most preferred robust solution.  
Moreover, the decision-maker can iteratively interact with the system, freely investigate different aspects of the problem, deepen their insight into the problem's complexity, and avoid significant losses due to the vast possibilities in the unknown future. 
Other unique features of the proposed decision support approach are reducing the cognitive load, structural learning, and decision-making environment to compare the consequences of the potential candidate solutions from different criteria and future realizations. 

Unfortunately, most predictions are not precisely realized in reality, and planning based on a simple average or the most probable scenario (the basis of probabilistic methods) is no longer enough in our fast-growing and continuously changing world and may fail because of a different realization. 
The primary purpose of this study is to support forest planners in better preparing for the unknown future by identifying the most sustainable and robust harvest schedule that moderately works well in a broad range of unknown futures. 

%limitations
% future research direction
As mentioned, this paper sheds light on dealing with the existing deep uncertainty in tactical forest harvest scheduling and providing better support for sustainable and robust planning. Nevertheless, many other aspects remain to be addressed in the future. 
For example, as pointed out in our experiments, although studying the whole planning horizon at once is required, the complexity it brings is high, making decision-making challenging. 
Furthermore, the robustness of the schedule is more important in the earlier periods than the later ones, as one can update the plan later as one moves on the planning horizon. 
One way is to extend the proposed approach to multi-stage decision-making problems, considering adaptable and contingency plans that can increase robustness \citep{Shavazipour2019,shavazipour2023}. 

% different sources of uncertainty like demand, accessibility of roads, etc
%The proposed decision support can also be adjusted for similar decision-making problems with different objective functions and/or sources of uncertainty. For instance, as discussed, having evenly distributed stands to be harvested in the schedule can be considered as an additional objective or constraint to avoid scheduling very few stands for harvesting if harvesting is limited for any reason, like route accessibility. 
%On the other hand, route accessibility can even be considered a deeply uncertain parameter in the model to be able to identify more sustainable and robust plans, taking into account this vital practical aspect in planning as well. 
%Combining route accessibility and multi-stage planning will draw an interesting future research direction for extending the current study. 

The proposed decision support can also be adjusted for similar decision-making problems with different objective functions and/or sources of uncertainty. For instance, as discussed, having evenly distributed stands to be harvested in the schedule can be considered as an additional objective or constraint to avoid scheduling very few stands for harvesting if harvesting is limited for any reason, like route accessibility. 
On the other hand, route accessibility can even be considered a deeply uncertain parameter in the model to be able to identify more sustainable and robust plans, taking into account this vital practical aspect in planning as well. 
Combining route accessibility and multi-stage planning will draw an interesting future research direction for extending the current study. This will also pave the way for considering demand uncertainty, leading to promising decision support for practical use, although further testing and validations are necessary.

\textcolor{blue}{
\bibliographystyle{agsm}
\bibliography{paper3refs}
}

\end{document}